\documentclass{amsart}
\usepackage{amssymb,amsmath,epsf}

\include{graphics}

\newtheorem{theorem}{Theorem}
\newtheorem{corollary}[theorem]{Corollary}
\newtheorem{definition}[theorem]{Definition}

\newtheorem{lemma}[theorem]{Lemma}
\newtheorem{remark}[theorem]{Remark}

\newcommand{\ZZ}{{\mathbb{Z}}}

\newcommand{\abs}[1]{\lvert#1\rvert}

\begin{document}

\baselineskip=18pt
\begin{center}
{\textbf{SOME IDENTITIES AND FORMULAS INVOLVING GENERALIZED
CATALAN NUMBERS}} \footnote{Supported by Norway-SA grant 2067063.

Address: School of Mathematical Sciences, University of
KwaZulu-Natal, Pietermaritzburg, 3209 South Africa \quad e-mail:
ngs@ukzn.ac.za }
\end{center}

\begin{center}

Siu-Ah Ng

\medskip

%\rm{December 2004}

\end{center}

\bigskip

\begin{quote}
\small{\em A generalization of the Catalan numbers is considered.
New results include binomial identities, recursive relations and a
close formula for the multivariate generating function. A simple
expression for the Catalan determinant is derived. }

\bigskip

\end{quote}

\vskip 30pt

\section{Introduction}\label{background}

\noindent By a \emph{path}, we mean a finite sequence
$\{a_i\}_{0\leq i\leq n}$ where $a_i\in\ZZ^2$ and either
$a_{i+1}=a_i+(1,0)$ or $a_{i+1}=a_i+(0,1).$ That is to say, a walk
from $a_0$ to $a_n$ whose steps consist of horizontal or vertical
positive unit movements.

The classical Catalan number is the number of paths from the
origin to $(n,n)$ without crossing the diagonal. It is given by
$\,\displaystyle{\frac{1}{n+1}\,{2n \choose n}}\,$ and has many
well-known interpretations (Cf \cite{HP1}). In \cite{Ng} we came
across the following natural generalization:

\bigskip

\begin{definition}\label{catalandefinition} Let $\,0\leq m\leq n.\,$
The generalized Catalan number is defined as the number of paths
from $(0,-2m)$ to $(n-m, n-m)$ without crossing the line $y=x$ and
is denoted by $\,\mathcal{C}_{n,m}.\,$
\end{definition}

\bigskip

For completeness, we give in \S \ref{catalang} a proof of
$\,\displaystyle{\mathcal{C}_{n,m}\,=\,\frac{2m+1}{n+m+1}\,{2n\choose
n+m}}.\,$ The number $\,\mathcal{C}_{n,m}\,$ has an ancient
history. In \cite{B} Bertrand studied $\,a_{n,m},\,$ the number of
paths from the origin to $(n, m)$ without crossing the diagonal,
as the number of solutions to the ballot problem. One verifies
that $\displaystyle{\,\mathcal{C}_{n,m}\,=\,a_{n+m,n-m},\quad
n\geq m.\,}$ (For another derivation, one can check that
$\displaystyle{\,\mathcal{C}_{n,m}\,=\,d_{(1-2m)(1+n-m)}}$ for the
$d_{qk}$ defined in \cite{HP2}.)

We then move on to proving new identities and formulas.

In \S \ref{variations} we consider some variations of
$\,\mathcal{C}_{n,m}\,$ and compute closed formulas for them. In
\S \ref{identities} we derive some binomial identities by counting
Catalan numbers and their variants. In \S \ref{recur} we obtain
recursive relations for the generalized Catalan numbers, from
which identities for the central binomial coefficients and an
identity for certain sum of products of the generalized Catalan
numbers are given. In the main section \S \ref{genfunction} we
give a close formula for the multivariate generating function for
the generalized Catalan numbers. In \S \ref{moreformulas} more
formulas are produced, including certain product of the generating
functions. We end in the final section \S \ref{det} with an
elegant formula for the \emph{Catalan determinant}.

\vskip 30pt

\section{Generalized Catalan numbers}\label{catalang}

For completeness, we give a direct derivation of a close formula
for $\,\mathcal{C}_{n,m}\,$ although as remarked earlier it is
already included in references such as \cite{B} and \cite{HP2}.

\begin{theorem}\label{catalan}
Given integers $\,0\leq m\leq n,\,$ we have
\begin{equation}\label{catalangeq}
\mathcal{C}_{n,m}\,=\,\frac{2m+1}{n+m+1}\,{2n\choose n+m}.
\end{equation}
\end{theorem}

\begin{proof} Fix $\,0\leq m\leq n.$
We let $\,\Lambda\,$ be the set of paths from $\,(0,-2m)\,$ to
\mbox{$\,(n-m,n-m)\,$} and $\,\Gamma\,$ be the subset of ``bad"
paths, i.e. those that cross the line $\,y=x.\,$ So
$\,\mathcal{C}_{n,m}\,=\,\abs{\Lambda}\,-\,\abs{\Gamma}.\,$

\begin{center}
\begin{figure}[h]
\hskip 150pt \epsfbox[265 120 50 10]{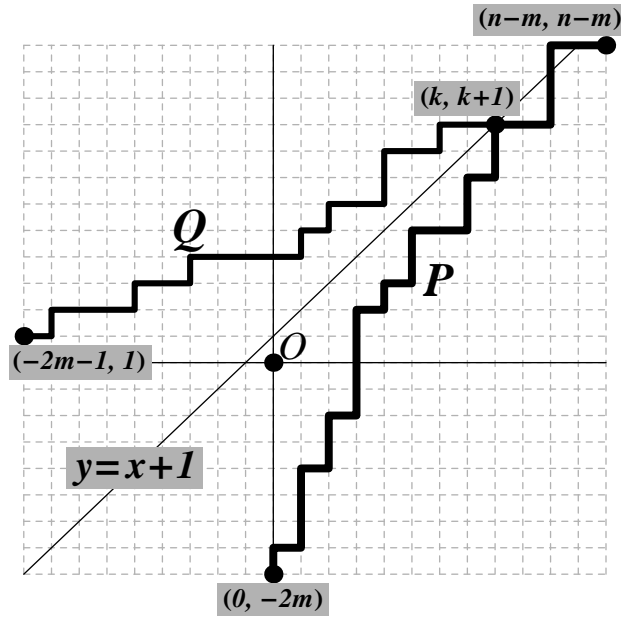}

\vskip 210pt

\hskip -145pt\caption{\label{fig1} Portion of $P$ reflected across
$y=x+1$ to form $Q.$}
\end{figure}
\end{center}
Given a path $\,P\in\Gamma,\,$ let $k$ be smallest such that $P$
crosses $\,y=x\,$ at $\,(k,k),\,$  i.e. both $(k,k)$ and $(k,k+1)$
belong to $P.$ Let $Q$ be the path obtained from $P$ by reflecting
across the line $y=x+1$ only the portion of $P$ from $\,(0,-2m)\,$
to $\,(k,k+1).\,$ The result is a path from $\,(-2m-1,1)\,$ to
\mbox{$\,(n-m,n-m).\,$} See Figure \ref{fig1} above. This is
simply the Andr\'e's reflection method and one can easily see that
it establishes a bijection between $\Gamma$ and the set of paths
from $\,(-2m-1,1)\,$ to \mbox{$\,(n-m,n-m).\,$} Therefore
\[\abs{\Gamma}\,=\,{2(n-m)-(-2m-1+1)\choose (n-m)-1}\,=\,{2n\choose n-m-1}.\]

On the other hand,
\[\abs{\Lambda}\,=\,{2(n-m)-(-2m)\choose n-m}\,=\,{2n\choose n-m}\,=\,{2n\choose n+m}.\]
Hence
\[\mathcal{C}_{n,m}\,=\,{2n\choose n+m}\,-\,{2n\choose n-m-1}\]
and the conclusion follows by noticing that
\[{2n\choose n-m-1}\,=\,\frac{n-m}{n+m+1}\,{2n\choose n-m}\,=\,\frac{n-m}{n+m+1}\,{2n\choose n+m}.\]

\end{proof}

\bigskip

\bigskip

\noindent \textbf{Bertrand's ballot problem}: Suppose the final
result of an election is such that party $A$ gets $2m$ more votes
than its opposition $B,$ say $A$ gets $n+m$ and $B$ gets $n-m$
votes. What is the probability that $A$ is at least $2m$ votes
ahead of $B$ in all partial counts throughout the whole process?

\bigskip

\begin{corollary}\label{unfairballot}
The probability of the event in the generalized ballot problem is
$\,\displaystyle{\frac{2m+1}{n+m+1}}.\,$
\end{corollary}

\begin{proof}
By shifting all paths $2m$ units up, it is clear that
$\mathcal{C}_{n,m}$ is the number of paths from $(0,0)$ to
$(n-m,n+m)$ without crossing $y=x+2m.$ Then identify each vote for
party $A$ with $(0,1)$ and that for party $B$ with $(1,0).$

The conclusion follows by noticing that each possible ballot
outcome in the problem corresponds to a path from $(0,0)$ to
$(n-m,n+m),$ hence there are $\,\displaystyle{{2n\choose n+m}}\,$
of them in total, and those ballot outcome with $A$ leading $B$ by
at least $2m$ votes corresponds to paths from $(0,0)$ to
$(n-m,n+m)$ without crossing $y=x+2m.$
\end{proof}

\vskip 60pt

\section{Some variations}\label{variations}
Superficially, the definition of $\,\mathcal{C}_{n,m}\,$ seems to
be too restricted, since it considers only the even number $2m$
and allowing touching $y=x.$ Here in this section we consider some
variations and they can be computed in similar ways.

\bigskip

\begin{definition} Let $\,0\leq m\leq n.\,$

$\bar{\mathcal{C}}_{n,m}\,:=\,$the number of paths from
$(0,-2m-1)$ to $(n-m,n-m)$ without crossing $y=x.$

$\mathcal{D}_{n,m}\,:=\,$the number of paths from $(0,-2m)$ to
$(n-m,n-m)$ strictly below $y=x$ until reaching $(n-m,n-m).$

$\bar{\mathcal{D}}_{n,m}\,:=\,$the number of paths from
$(0,-2m-1)$ to $(n-m,n-m)$ strictly below $y=x$ until reaching
$(n-m,n-m).$
\end{definition}

\bigskip

\begin{theorem}\label{catalanvariations}
Suppose $\,0\leq m\leq n.\,$ Then
\begin{enumerate}
    \item $\displaystyle{\bar{\mathcal{C}}_{n,m} \,=\,\mathcal{D}_{n+1,m+1}
    \,=\, \frac{2m+2}{n+m+2}\,{2n+1\choose n-m};}$
    \medskip
    \item $\displaystyle{\mathcal{D}_{n,m} \,=\, \frac{2m}{n+m}\,{2n-1\choose n-m};}$
    \medskip
    \item $\displaystyle{\bar{\mathcal{D}}_{n,m} \,=\,\mathcal{C}_{n,m}\,=\, \frac{2m+1}{n+m+1}\,{2n\choose n-m}.}$
\end{enumerate}
\end{theorem}

\begin{proof}
We first prove (2). Notice that a path $P$ from the set defining
$\mathcal{D}_{n,m}$ must pass $(n-m,n-m-1)$ before reaching
$(n-m,n-m).$ Similar to the proof of Theorem \ref{catalan}, $P$
can be reflected across the line $y=x$ to form a path $Q$ from
$\,(-2m,0)\,$ to \mbox{$\,(n-m,n-m-1),\,$} as shown in Figure
\ref{fig2}.

\newpage

\begin{center}
\begin{figure}[h]
\hskip 150pt \epsfbox[280 120 50 10]{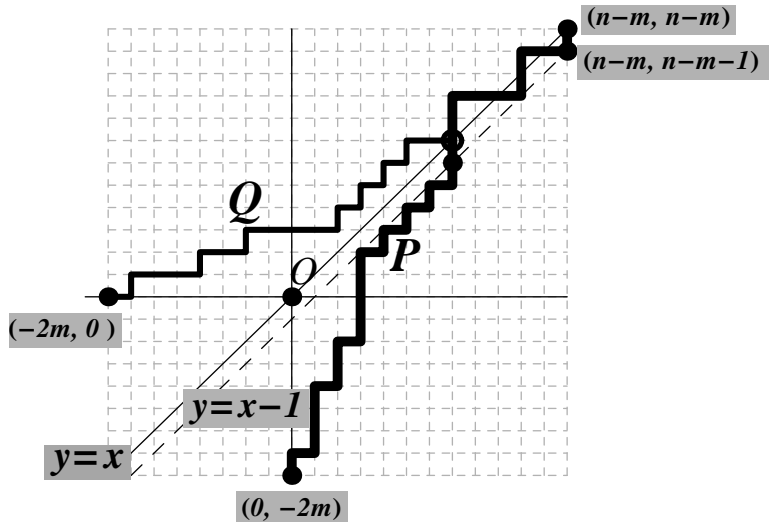}

\vskip 210pt

\hskip -145pt\caption{\label{fig2} Reflection across $y=x.$}
\end{figure}
\end{center}
We let $\,\Lambda\,$ be the set of paths from $\,(0,-2m)\,$ to
\mbox{$\,(n-m,n-m-1)\,$} and $\,\Gamma\,$ be the subset of ``bad"
paths, i.e. those that cross the line $\,y=x-1.\,$ Clearly,
$\,\displaystyle{\abs{\Lambda}\,=\,{2n-1\choose n-m}.}\,$ By the
reflection trick used in the proof of \mbox{Theorem
\ref{catalan},} we can identify paths in $\,\Gamma\,$ with paths
from $(-2m,0)$ to $(n-m,n-m-1).$ Therefore
$\,\displaystyle{\abs{\Gamma}\,=\,{2n-1\choose n+m}.}\,$

So
$\,\displaystyle{\mathcal{D}_{n,m}\,=\,\abs{\Lambda}\,-\,\abs{\Gamma}\,=\,{2n-1\choose
n-m}-{2n-1\choose n+m}\,=\,\frac{2m}{n+m}\,{2n-1\choose n-m}.}\,$

\bigskip

To prove (1), by shifting $1$ unit down, we can identify paths
defining $\,\bar{\mathcal{C}}_{n,m}\,$ with paths from $(0,-2m-2)$
to $(n-m,n-m)$ strictly below $y=x$ until reaching $(n-m,n-m).$
Therefore
\[\bar{\mathcal{C}}_{n,m}\,=\,\mathcal{D}_{n+1,m+1}\]
and the conclusion follows.

Now for (3), we notice that by shifting $1$ unit up, paths
defining $\,\bar{\mathcal{D}}_{n,m}\,$ are identified with paths
from $(0,-2m)$ to $(n-m,n-m)$ without crossing $y=x.$ Therefore
$\;\bar{\mathcal{D}}_{n,m}\,=\,\mathcal{C}_{n,m}.$
\end{proof}

\bigskip

\section{Identities}\label{identities}

In this section we derive some interesting identities which are
corollaries of \mbox{Theorem \ref{catalan}} or Theorem
\ref{catalanvariations}.

\bigskip

\begin{corollary}\label{id1}
\begin{equation}\label{identityeq1}
\sum_{k=m}^{n-1}\,\frac{1}{2k+1}\,{2k+1\choose
k-m}\,{2(n-k)\choose n-k}\;=\;\frac{2}{2m+1}{2n\choose n-m-1},
\end{equation}
where $\,0\,\leq\,m\,\leq n-1.$
\end{corollary}

\begin{proof}
A ``bad" path $P$ from $(0,-2m)$ to \mbox{$(n-m, n-m)$}, i.e. one
that crosses $y=x,$ is a path from $(0,-2m)$ to some $(k, k),$
immediately before crossing $y=x,$ then proceed to $(k,k+1)$ and
takes any path to $(n-m, n-m)$ as shown in Figure~\ref{fig3}
below.

\begin{center}
\begin{figure}[h]
\hskip 150pt \epsfbox[280 150 50 10]{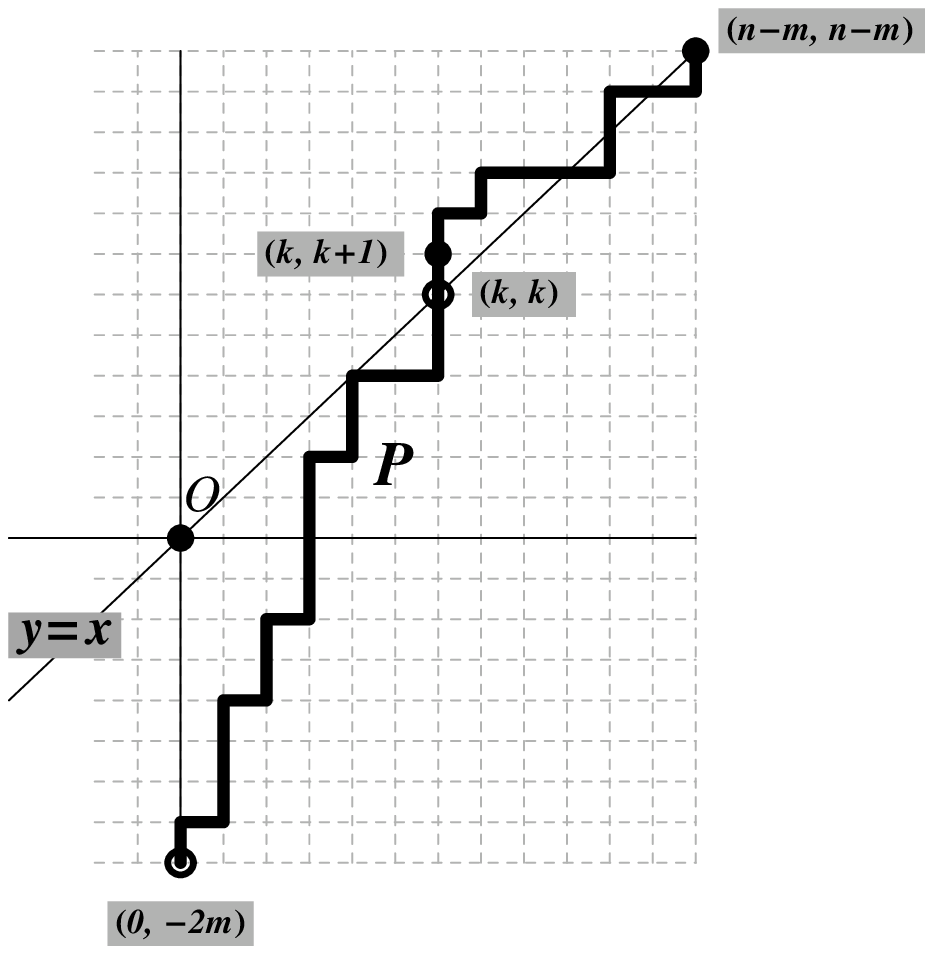}

\vskip 270pt

\hskip -145pt\caption{\label{fig3} Decomposition of a ``bad" path
$P.$}
\end{figure}
\end{center}

The number of paths from $(0,-2m)$ to \mbox{$(k,
k)\,=\,((m+k)-m,(m+k)-m)$} without crossing $y=x$ is given by
$\,\mathcal{C}_{m+k,m},\,$ while the number of paths from
$(k,k+1)$ to $(n-m, n-m)$ is given by
$\,\displaystyle{{2(n-m-k)-1\choose n-m-k}}.\,$ Hence the number
of the ``bad" paths are

\bigskip

\begin{align*}
&\quad\, \sum_{k=0}^{n-m-1}\mathcal{C}_{m+k,m}\,{2(n-m-k)-1\choose
n-m-k}\\
&=\,\sum_{k=0}^{n-1}\mathcal{C}_{k,m}\,{2(n-k)-1\choose
n-k}\\
&=\,\sum_{k=0}^{n-1}\,\frac{2m+1}{k+m+1}\,{2k\choose
k-m}\,{2(n-k)-1\choose n-k}.
\end{align*}

\bigskip

Since this number is also $\displaystyle{{2n\choose n-m}
-\mathcal{C}_{n,m}}$ the conclusion follows from the last part of
the proof of Theorem \ref{catalan}.
\end{proof}

\bigskip

\bigskip

\begin{corollary}\label{id3}
\begin{equation}\label{identityeq3}
\sum_{k=m}^n\,\frac{1}{k(n-k+1)}\,{2k\choose k-m}\,{2(n-k)\choose
n-k}\,=\,\frac{2m+1}{m(n+m+1)}\,{2n\choose n-m}
\end{equation}
where $\,1\,\leq\,m\,\leq n-1.$
\end{corollary}

\begin{proof}
Note that a path from $(0,-2m)$ to $(n-m, n-m)$ without crossing
the line $y=x$ can be decomposed into a path from $(0,-2m)$ to the
first touch at $y=x,$ some $(k, k),$ $\,0\leq k\leq n-m,$ followed
by a path from $(k, k),$ to $(n-m, n-m)$ without crossing the line
$y=x,$ as shown in Figure~\ref{fig4} below.

\newpage

\begin{center}
\begin{figure}[h]
\hskip 150pt \epsfbox[260 120 50 10]{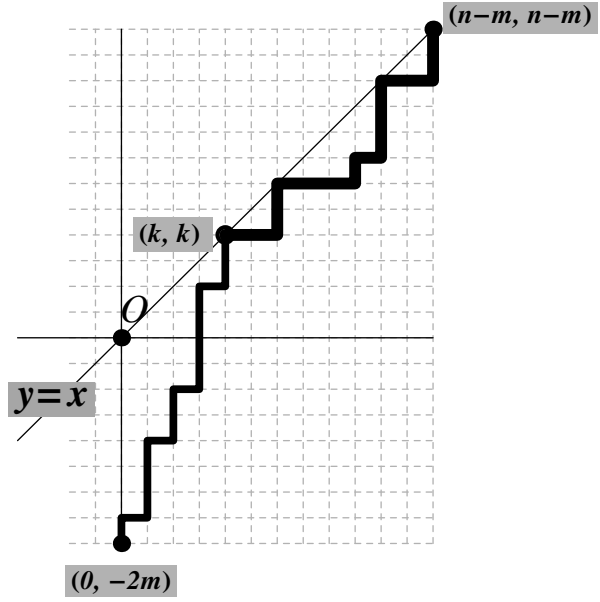}

\vskip 210pt

\hskip -80pt\caption{\label{fig4} The decomposition at $(k,k).$}
\vskip -20pt
\end{figure}
\end{center}
The number of the former paths is $\mathcal{D}_{k+m,m}$ and the
number for the latter is $\mathcal{C}_{n-m-k,0},\,$ therefore we
have
\begin{equation}\label{identityeq3a}
\mathcal{C}_{n,m}\,=\,\sum_{k=0}^{n-m}\,\mathcal{D}_{k+m,m}\,\mathcal{C}_{n-m-k,0}.\end{equation}
Hence, by Theorem \ref{catalan} and Theorem
\ref{catalanvariations},
\begin{eqnarray*}
  \frac{2m+1}{n+m+1}\,{2n\choose n-m} &=& \sum_{k=0}^{n-m}\frac{2m}{k+2m}\,
  {2(k+m)-1\choose k}\,\frac{1}{n-m-k+1}{2(n-m-k)\choose n-m-k}\, \\
  \, &=& \sum_{k=m}^{n}\frac{2m}{k+m}\,{2k-1\choose k-m}\,\frac{1}{n-k+1}\,{2(n-k)\choose n-k} \\
     &=& \sum_{k=m}^{n}\frac{m}{k(n-k+1)}\,{2k\choose k-m}\,{2(n-k)\choose
     n-k},
\end{eqnarray*}
and equation (\ref{identityeq3}) follows.
\end{proof}

\bigskip

\begin{remark}
It is easy to check that equation (\ref{identityeq1}) can be
transformed into the following form:
\begin{multline}\label{identityeq11}
    \sum_{k=m}^{n}\,\frac{2n-2k+1}{(k+m+1)(n-k+1)}\,{2k\choose
    k-m}\,{2(n-k)\choose n-k}\\
    =\;\frac{(2n+1)(2n+2)}{(2m+1)(n+m+1)(n+m+2)}{2n\choose n-m}.
\end{multline}
Although equations (\ref{identityeq3}) and (\ref{identityeq11})
are similar, it is not clear how they are related.
\end{remark}

\vfill

\begin{corollary}\label{id2}
\begin{equation}\label{identityeq2}
\sum_{h=0}^{n-m}\,\frac{k}{h+k}\,{2h+k-1\choose
h}\,{2n-2h-k\choose n-m-h}\,=\,{2n\choose n-m}
\end{equation}
where $\,0\,\leq\,m\,\leq n-1\,$ and $\,1\leq k\leq 2m.$
\end{corollary}

\begin{proof}
The number of paths from $(0,-2m)$ to $(n-m,n-m)$ is
$\,\displaystyle{{2n\choose n-m}},\,$ the right side of the above
equation (\ref{identityeq2}).

For each $\,1\leq k\leq 2m\,$ we can decompose these paths into
the path from $(0,-2m)$ first touching the line $y=x-2m+k$ at
$(h,h-2m+k),$ for some $0\leq h\leq n-m,$ followed by a path from
$(h,h-2m+k)$ to $(n-m,n-m),$ as shown in Figure~\ref{fig5} below.

\begin{center}
\begin{figure}[h]
\hskip 150pt \epsfbox[260 100 50 10]{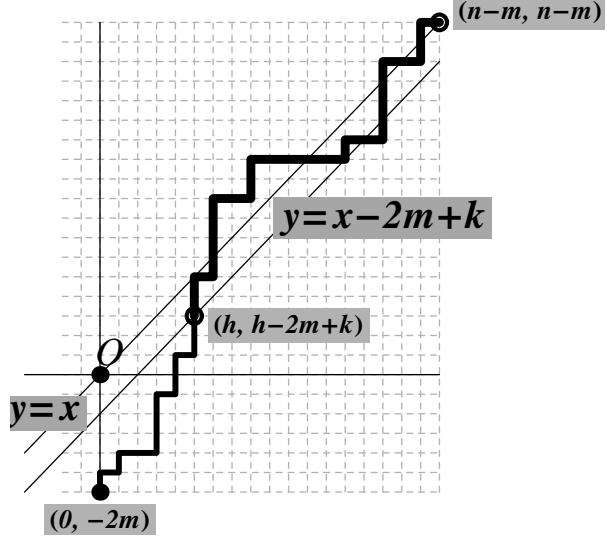}

\vskip 180pt

\hskip -80pt\caption{\label{fig5} The decomposition at
$(h,h-2m+k).$}
\end{figure}
\end{center}

The number of paths from $(h,h-2m+k)$ to $(n-m,n-m)$ is
$\,\displaystyle{{2n-2h-k\choose n-m-h}}.\,$

By shifting up $\,2m-k\,$ units, we can see that a path from
$(0,-2m)$ first touching the line $y=x-2m+k\,$ at $\,(h,h-2m+k)\,$
corresponds to a path from $(0,-k)$ first touching the line
$y=x\,$ at $\,(h,h).\,$ Therefore, if $k$ is even, the number is
given by
\[\mathcal{D}_{h+\frac{k}{2},\frac{k}{2}}\,=\,\frac{k}{h+k}\,{2h+k-1\choose h};\]
if $k$ is odd, the number is given by
\[\bar{\mathcal{D}}_{h+\frac{k-1}{2},\frac{k-1}{2}}\,=\,\frac{k}{h+k}\,{2h+k-1\choose h};\]
both producing the same number.

Summing the product $\displaystyle{\frac{k}{h+k}\,{2h+k-1\choose
h}\,{2n-2h-k\choose n-m-h}}$ over $0\leq h\leq n-m,$ equation
(\ref{identityeq2}) follows.
\end{proof}

\bigskip

\begin{remark}\label{remarkidentityeq2}
It appears that equation (\ref{identityeq2}) actually holds for
all $\,k>0.$ The current proof doesn't seem to work for $\,k>2m.$
\end{remark}

\vskip 50pt

\section{Recursive relations and identities}\label{recur}

\noindent In this section we will derive two recursive relations
for the generalized Catalan numbers. As corollaries, we obtain
some identities for the central binomial coefficients. The first
recursive relation for the generalized Catalan numbers we consider
is the following.

\begin{theorem}\label{reqthm2} For $\,0\leq m<n,$ we have the
following recursive relation
\begin{equation}\label{thm2eq}
{\mathcal{C}_{n,m}}\,=\,\frac{n(2m+1)}{(n-m)(m+1)}\,\sum_{k=0}^{n-m-1}\,{\mathcal{C}_{k+m,\,m}}
\,{\mathcal{C}_{n-m-k-1,\,0}}.
\end{equation}
\end{theorem}

\begin{proof}
A path defining $\,{\mathcal{C}_{n,m}}\,$ is either a path
touching the line $\,y=x\,$ only once at the terminal point
$\,(n-m,n-m)\,$ (the number of such is given by
$\,{\mathcal{D}_{n,m}}\,$ ) or a path touching some $\,(k,k),\,$
where \mbox{$\,0\leq k\leq n-m-1\,$} for the last time before
hitting the terminal point \mbox{$\,(n-m,n-m),\,$} as shown in
Figure~\ref{fig11}.
\begin{center}
\begin{figure}[h]
\hskip 150pt \epsfbox[280 120 50 10]{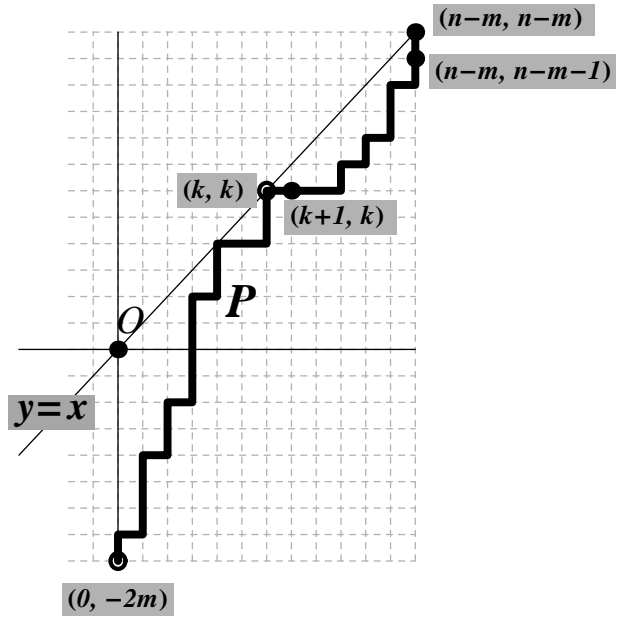}

\vskip 210pt

\hskip -145pt\caption{\label{fig11} A path $P$ touching $y=x$ at
$(k,k)$ for the last time before touching the terminal point
$\,(n-m,n-m).$}
\end{figure}
\end{center}
Write $\,(k,k)\,$ as $((k+m)-m,(k+m)-m),\,$ we have
\[{\mathcal{C}_{n,m}}\,=\,{\mathcal{D}_{n,m}}\,+\,\sum_{k=0}^{n-m-1}\,{\mathcal{C}_{k+m,\,m}}
\,{\mathcal{C}_{n-m-k-1,\,0}}.\] But
\[{\mathcal{D}_{n,m}}\,=\,\frac{2m}{n+m}\,{2n-1\choose n-m}\,=\,\frac{m(n+m+1)}{n(2m+1)}\,{\mathcal{C}_{n,m}},\]
hence the result follows.
\end{proof}

\bigskip

As a consequence, we have an identity for the $n$th central
binomial coefficient.

\begin{corollary}\label{cor1}
Let $\,0\leq m<n,$ then \begin{equation}\label{cor1eq} {2n\choose
n}\,=\,\frac{n+m+1}{m+1}\,\left(\prod_{k=0}^m\frac{n+k}{n-k}\right)\,\sum_{k=0}^{n-m-1}\,{\mathcal{C}_{k+m,\,m}}
\,{\mathcal{C}_{n-m-k-1,\,0}}.
\end{equation}
\end{corollary}

\begin{proof}
By equation (\ref{thm2eq}) in Theorem \ref{reqthm2} we have:
\[\frac{2m+1}{n+m+1}\,\frac{(2n)!}{(n+m)!\,(n-m)!} \,=\,
\frac{n(2m+1)}{(n-m)(m+1)}\,\sum_{k=0}^{n-m-1}\,
  {\mathcal{C}_{k+m,\,m}}
\,{\mathcal{C}_{n-m-k-1,\,0}}, \] i.e.
\[\frac{(2n)!}{(n+m)!\,(n-m)!} \,=\,
\frac{n+m+1}{m+1}\,\frac{n}{n-m}\,\sum_{k=0}^{n-m-1}\,
   {\mathcal{C}_{k+m,\,m}}
\,{\mathcal{C}_{n-m-k-1,\,0}},\] hence the result follows if
$m=0;$ for $m>0,$ we get
\[\frac{(2n)!}{n!^2} \,=\, \frac{n+m+1}{m+1}\,\frac{n(n+1)\cdots(n+m)}{n(n-1)\cdots(n-m+1)(n-m)}\,\sum_{k=0}^{n-m-1}\,
   {\mathcal{C}_{k+m,\,m}}
\,{\mathcal{C}_{n-m-k-1,\,0}},\] and equation (\ref{cor1eq}) is
proved.
\end{proof}

\bigskip

The following is an identity to be used in \S \ref{moreformulas}.

\bigskip

\begin{corollary}\label{cor1a}
Let $\,1\leq m<n,$ then
\begin{eqnarray}\label{cor1aeq} \sum_{k=0}^{n-m-1}\,{\mathcal{C}_{k+m,\,m}}
\displaystyle{{\mathcal{C}_{n-m-k-1,\,0}}}\,&=&\,\displaystyle{\frac{m+1}{n!}\,\left(\prod_{k=1}^m\,(n-k)\right)\,
\left(\prod_{k=m+2}^n\,(n+k)\right)}\\
&=&\,\displaystyle{\left(\prod_{k=1}^m\,\frac{n-k}{k}\right)\,
\left(\prod_{k=m+2}^n\,\frac{n+k}{k}\right).}\notag
\end{eqnarray}
\end{corollary}

\begin{proof}
This is just another form of equation (\ref{cor1eq}) in Corollary
\ref{cor1} by noticing that
\[{2n\choose n}\,=\,\prod_{k=1}^n\,\frac{n+k}{k}.\]
\end{proof}

\bigskip

Now we consider another recursive relation for the generalized
Catalan numbers.

\bigskip

\begin{theorem}\label{recursive}
Let $\,0\leq m\leq n.$ Then
\begin{equation}\label{thm1eq}
{\mathcal{C}_{n,m}}\,=\,\frac{m}{2m+1}\,\sum_{k=m}^n\,\frac{k+m+1}{k}\,
{\mathcal{C}_{k,m}}\,{\mathcal{C}_{n-k,0}}.
\end{equation}
\end{theorem}

\begin{proof}
From Theorem 5 in \cite{Ng} we can obtain
\[\mathcal{D}_{n,m}\,=\,\frac{m(n+m+1)}{n(2m+1)}\,{\mathcal{C}_{n,m}}.\]
The result then follows from the identity (4) of Corollary 7 in
\cite{Ng}.
\end{proof}

From the above we have another expression for the $n$th central
binomial coefficient. The proof is similar to that of Corollary
\ref{cor1}.

\begin{corollary}\label{cor6}
Let $\,1\leq m\leq n,$ then
\begin{equation}\label{cor6eq} {2n\choose
n}\,=\,\frac{m(n-m)}{(2m+1)^2}\,\left(\prod_{k=1}^{m+1}\frac{n+k}{n-k+1}\right)\,
\sum_{k=m}^{n}\,\frac{k+m+1}{k}\,{\mathcal{C}_{k,\,m}}\,{\mathcal{C}_{n-k,\,0}}.
\end{equation}\qed
\end{corollary}

\vskip 30pt

\section{Generating functions}\label{genfunction}

\noindent In this section we derive a formula for the generating
function of generalized Catalan numbers. A comprehensive teatment
of generating functions can be found in \cite{FS}.

For fixed $\,m\geq 0\,$ we let $\,\gamma_m (x)\,$ denote the
generating function of $\,{\mathcal{C}_{n+m,\,m}},\,$ i.e. the
power series
\begin{equation}\label{gendef}
 \gamma_m (x)\,=\,\sum_{n=0}^\infty\,{\mathcal{C}_{n+m,\,m}}\,x^n.
\end{equation}
It is well-known for the classical Catalan numbers
$\,{\mathcal{C}_{n,\,0}}\,$ that
\[\gamma_0 (x)\,=\,\frac{1-\sqrt{1-4x}}{2x}.\]
(See for example \cite{HP1}.)

We denote the corresponding multivariate generating function by:
\begin{equation}\label{xygendef}
 \Gamma (x, y)\,=\,\sum_{n, m=0}^\infty\,{\mathcal{C}_{n+m,\,m}}\,x^n\,y^m.
\end{equation}

In \cite{CR} a formula for the generating function for
$\,a_{n,m}\,$ is given, \emph{i.e.}
\[\,a(x,y)\,:=\,\sum_{n\geq m\geq 0}a_{n,m} x^n
y^m\,=\,\frac{1-y\,\gamma_0 (xy)}{1-x-y}.\] Using the relation
$\,\mathcal{C}_{n,m}\,=\,a_{n+m,n-m},\,$ it can be verified that
\[\Gamma
(x,y)\,=\,\frac12\Big(a\Big(\sqrt{y},\frac{x}{\sqrt{y}}\Big)\,+\,a\Big(-\sqrt{y},
-\frac{x}{\sqrt{y}}\Big)\Big).\] However we have another
derivation here:

\begin{theorem}\label{mainthm}
\begin{equation}\label{xygenthm}
\Gamma (x,y)\,=\,\frac{\gamma_0 (x)}{1\,-\,y\, (\gamma_0
(x))^2}\quad\text{i.e.}\quad
\frac{x(1-\sqrt{1-4x}\,)}{2x(x+y)-y(1-\sqrt{1-4x}\,)}.\end{equation}\qed
\end{theorem}

\bigskip

We first need the following lemmas:

\begin{lemma}\label{mainlem}
Let $n\geq 0\,$ then
\begin{equation}\label{mainlemeq}
{\mathcal{C}_{n+m+1,\,m+1}}\,=\,\sum_{a+b+c\,=\,n}\,{\mathcal{C}_{a+m,\,m}}\,{\mathcal{C}_{b,\,0}}\,{\mathcal{C}_{c,\,0}},
\end{equation}
where we sum over $\,a,\,b,\,c\,\geq 0.$\end{lemma}

\begin{proof}
We count paths defining ${\mathcal{C}_{n+m+1,\,m+1}},$ i.e. paths
from $(0, -2(m+1))$ to $(n,n)=(n+m+1-(m+1),n+m+1-(m+1))$ without
crossing $y=x,$ by using their first touch $(a, a-1)$ and
$(a+b,a+b)$ at the lines $y=x-1$ and $y=x$ respectively as in
Figure~\ref{fig22} below.

\vskip 50pt

\begin{center}
\begin{figure}[h]
\hskip 150pt \epsfbox[280 120 50 10]{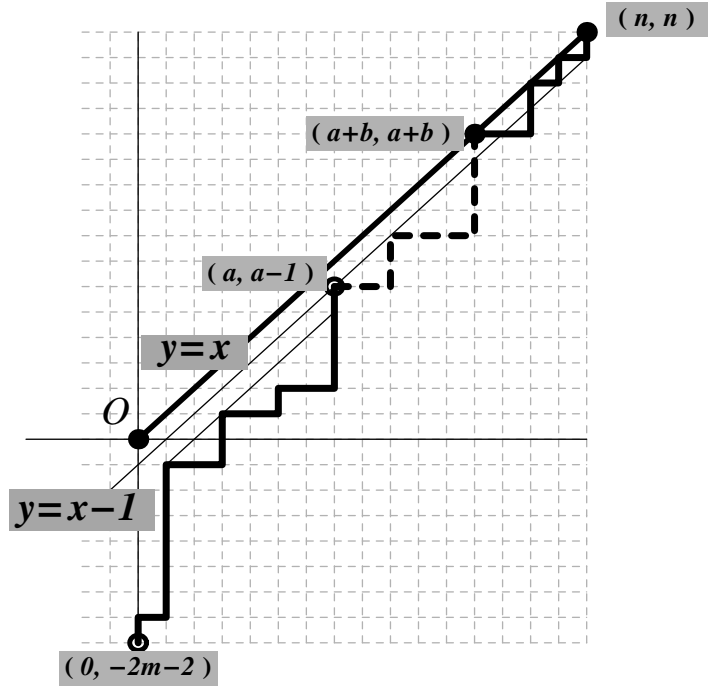}

\vskip 210pt

\hskip -135pt\caption{\label{fig22} A path defining
${\mathcal{C}_{n+m+1,\,m+1}}$ whose first touch at $y=x-1$ is
$(a,a-1)$ and whose first touch at $y=x$ is $(a+b,a+b)$ for some
$a,b\geq 0\,$ and $a+b\leq n.$ }
\end{figure}
\end{center}

The number of paths from $(0,-2m-2)$ that first touch $y=x-1$ at
$(a, a-1)$ is the same as those reaching $(a, a-2)$ without
crossing $y=x-2,$ i.e. the number is ${\mathcal{C}_{a+m,\,m}.}$
Similarly, the number of paths from $(a,a-1)$ first reaching $y=x$
at $(a+b,a+b)$ is ${\mathcal{C}_{b,0}}$  and the number of paths
from $(a+b,a+b)$ to $(n,n)$ without crossing $y=x$ is
${\mathcal{C}_{c,0}},$ where $c=n-a-b,$ hence the result is
proved.
\end{proof}

\bigskip

Now we have a recursive relation for the generating function
$\gamma_m (x):$

\begin{lemma}\label{gammalem}
For $m\geq 0,\,$ we have
\begin{equation}\label{gammalemeq}
\gamma_{m+1} (x)\,=\,\gamma_m (x)\,(\gamma_0 (x))^2.
\end{equation}
\end{lemma}

\begin{proof}
We use Lemma \ref{mainlem}:
\begin{eqnarray*}
  \gamma_{m+1} (x) &=& \sum_{n=0}^\infty\,{\mathcal{C}_{n+m+1,\,m+1}}\,x^n \\
   &=& \sum_{n=0}^\infty\, \sum_{a+b+c=\,n}\,{\mathcal{C}_{a+m,\,m}}\,
   {\mathcal{C}_{b,\,0}}\,{\mathcal{C}_{c,\,0}}\,x^n\\
   &=& \left(\sum_{n=0}^\infty\,{\mathcal{C}_{n+m,\,m}}\,x^n\right)\,
   \left(\sum_{n=0}^\infty\,{\mathcal{C}_{n,\,0}}\,x^n\right)^2\\
   &=& \gamma_m (x)\,(\gamma_0 (x))^2.
\end{eqnarray*}
\end{proof}

\bigskip

Now we prove Theorem \ref{mainthm} from Lemma \ref{gammalem}. Note
that:
\begin{eqnarray*}
  y\, \Gamma (x,y)(\gamma_0 (x))^2 &=& \sum_{m=0}^\infty\,\gamma_m (x)\,(\gamma_0 (x))^2\,y^{m+1} \\
   &=& \sum_{m=0}^\infty\,\gamma_{m+1} (x)\,y^{m+1}  \\
   &=& \Gamma (x,y)\,-\,\gamma_0 (x),
\end{eqnarray*}
from which the result follows.

\vskip 50pt

\section{More formulas}\label{moreformulas}
\noindent In this section we prove a few more fomulas relating to
the generating functions.

\bigskip

\begin{theorem}\label{catgenthm}
Let $\,m\geq 1.\,$ The generating function $\,\gamma_m (x)\,$ is
equal to

\[\frac{2x}{1-\sqrt{1-4x}}\,\left(1+\sum_{n=1}^\infty\,\left(\prod_{k=1}^m\frac{n+m-k+1}{k}\right)\,
\left(\prod_{k=m+2}^{n+m+1}\frac{n+m+k+1}{k}\right)\,x^n\right).\]

\end{theorem}

\medskip

\begin{proof}
Replacing $\,n-m-1\,$ by $\,n\,$ in Corollary \ref{cor1a} equation
(\ref{cor1aeq}), we obtain for $\,n\geq 1\,$ that
\[\sum_{k=0}^n\,{\mathcal{C}_{k+m,\,m}}\,{\mathcal{C}_{n-k,\,0}}\,=\,\left(\prod_{k=1}^m\,\frac{n+m-k+1}{k}\right)\,
\left(\prod_{k=m+2}^{n+m+1}\,\frac{n+m+k+1}{k}\right).\] Then the
theorem follows from
\begin{eqnarray*}
  \gamma_m (x)\,\gamma_0 (x) &=& \sum_{n=0}^\infty\,\sum_{k+h=n}\,{\mathcal{C}_{k+m,\,m}}\,{\mathcal{C}_{h,\,0}}\,x^n \\
   &=& \sum_{n=0}^\infty\,\left(\sum_{k=0}^n\,{\mathcal{C}_{k+m,\,m}}\,{\mathcal{C}_{n-k,\,0}}\right)\,x^n \\
   &=& 1\,+\,\sum_{n=1}^\infty\,\left(\prod_{k=1}^m\frac{n+m-k+1}{k}\right)\,
\left(\prod_{k=m+2}^{n+m+1}\frac{n+m+k+1}{k}\right)\,x^n.
\end{eqnarray*}
\end{proof}

\bigskip

We now the product $\,\gamma_m (x)\,\gamma_0 (x)\,$ and the
following power series:
\[\theta(x)\,=\,\sum_{n=0}^\infty\,\frac{(n+1)(m+1)}{(n+m+1)(2m+1)}\,{\mathcal{C}_{n+m+1,\,m}}\,x^n.\]

\bigskip

\begin{theorem}\label{thm2}
Let $\,0\leq m\leq n.\,$  Then
\[\theta(x)\,=\,\gamma_m (x)\,\gamma_0 (x)\,=\,\gamma_m (x)\,\frac{1-\sqrt{1-4x}}{2x}.\]
\end{theorem}

\begin{proof}
From Theorem \ref{reqthm2} equation (\ref{thm2eq}) we have:
\[{\mathcal{C}_{n+m+1,\,m}}\,=\,\frac{(n+m+1)(2m+1)}{(n+1)(m+1)}\,\sum_{k+h=n}\,{\mathcal{C}_{k+m,\,m}}\,
{\mathcal{C}_{h,\,0}},\] therefore
\[\theta(x)\,=\,\sum_{n=0}^\infty \left(\sum_{k+h=n}\,{\mathcal{C}_{k+m,\,m}}\,
{\mathcal{C}_{h,\,0}}\right)\,x^n\,=\;\gamma(x)\,\frac{1-\sqrt{1-4x}}{2x}.\]
\end{proof}

\bigskip

\begin{corollary}\label{anothercor}
For $\,1\leq m\leq n,\,$ we have
\begin{equation}\label{anothereq}
{\mathcal{C}_{n+m,\,m}}\,=\,\frac{(n+m)(2m+1)}{n(m+1)}\,\left(\prod_{k=1}^m\frac{n+m-k}{k}\right)\,
\left(\prod_{k=m+2}^{n+m}\frac{n+m+k}{k}\right).
\end{equation}
\end{corollary}

\begin{proof}
Compare the power series expansion of $\,\gamma_m (x)\,\gamma_0
(x)\,$ given in the last line of the proof of Theorem
\ref{catgenthm} and the one given by $\,\theta(x).\,$ Then
equation (\ref{anothereq}) is obtained by replacing $\,n+1\,$ for
$\,n.$
\end{proof}

\vskip 30pt

\section{A beautiful determinant}\label{det}

\noindent In this final section, we prove the following:

\begin{theorem}\label{thm1}
For integer $\,N\,\geq 0,$

\[\det\left[\frac{(2i+1)\,(2(i+j))!}{j!\,(2i+j+1)!}\right]_{0\leq i, j\leq N}\,=\,2^{\frac{N (N+1)}{2}}.\]
\qed
\end{theorem}

\bigskip

We suggest the name \emph{Catalan determinant} for the above, for
the entries of the matrix are the generalized Catalan numbers
$\,\mathcal{C}_{i+j,\, i}.\,$

First we need a slight modification of  \cite{W} Theorem 1. (See
also \cite{K}.) An almost identical proof is given here for
completeness.

\begin{lemma}\label{lem}
Let $\,f(x)\,=\,1+a_1 x+a_2 x^2+\dots\,$ be a formal power series
and let \mbox{$\,c_{i,j}\,=\,[ x^j ] f^{2i+1} (x),\,$} where
$\,[x^j]\,$ denotes the coefficient of $x^j$ in the series. Then
\[\det\left[c_{i,j}\right]_{0\leq i, j\leq N}\,=\,(2a_1)^{\frac{N (N+1)}{2}}.\]
\end{lemma}

\begin{proof}
Let
\[C\,:=\left[c_{i,j}\right]_{0\leq i, j\leq N}\quad\text{and}\quad\,B\,:=\left[b_{i,j}\right]_{0\leq i, j\leq N}\quad
\text{where}\quad b_{i,j}\,=\,(-1)^{i+j}{i\choose j}.\] The
$\,(i,k)$-entry of $\, BC\,$ is given by
\begin{eqnarray*}
  \sum_{j=0}^N b_{i,j}\,c_{j,k}&=& (-1)^i\,\,[ x^k ]
  \,\sum_{j=0}^N (-1)^j\,{i\choose j}\,f^{2j+1} (x)\\
   &=& ((-1)^i\,\,[ x^k ]
  \,\sum_{j=0}^i (-1)^j\,{i\choose j}\,f^{2j+1} (x)\quad\text{(since}\; {i\choose j}=0\;\text{for}\; j>i)\\
     &=& (-1)^i\,\,[ x^k ]
  \,\left(\left(1\,-\,f^2(x)\right)^i\,f(x)\right)\\
   &=& (-1)^i\,\,[ x^k ]
  \,\left(\left(-2a_1 x+\dots\right)^i\,(1+a_1 x+\dots)\right),
\end{eqnarray*}
which equals $\,(2a_1)^i\,$ if $\,k=i\,$ and equals $\,0\,$ if
$\,k<i,$ hence $\,BC\,$ is an upper triangular matrix with
diagonal entries $\,(2a_1)^i.\,$  In particular \[\,\det
BC\,=\,\prod_{i=0}^N\,(2a_1)^i\,=\,(2a_1)^{\frac{N (N+1)}{2}}.\]
But $\,B\,$ is lower triangular with diagonal entries $\,1,\,$ so
$\,\det B\,=\,1\,$ and the result follows.
\end{proof}

\bigskip

To prove Theorem \ref{thm1},  we let, as in \S \ref{genfunction},
$\,\gamma_i (x)\,$ be the generating function of
$\,{\mathcal{C}_{i+j,\,i}},\,$ i.e. the power series
\[ \gamma_i (x)\,=\,\sum_{n=0}^\infty\,{\mathcal{C}_{i+j,\,i}}\,x^j.
\]
Then
\[\gamma_0 (x)\,=\,\frac{1-\sqrt{1-4x}}{2x}\,=\,1+x+2x^2+5x^3+14x^4+\dots\] and Lemma~\ref{gammalem} gives that
\[ \gamma_i (x)\,=\,\gamma_0^{2i+1} (x)\quad\text{for}\quad i\geq
0.\] Therefore
\[{\mathcal{C}_{i+j,\,i}}\,=\,[ x^j ] \gamma_i (x)\,=\,[ x^j ] \gamma_0^{2i+1} (x)\]
and Theorem \ref{thm1} is now proved by applying Lemma \ref{lem}
to $\,f(x)\,=\,\gamma_0(x).$

\vskip 50pt

\bibliographystyle{amsplain}

\vfill

\end{document}